\documentclass[11pt]{amsart}
\usepackage{amssymb}
\usepackage{amsfonts}
\usepackage{amscd}
\usepackage{amsmath}
\usepackage{mathrsfs}
\usepackage{curves}
\usepackage{epsfig}
\usepackage[pass]{geometry}
\usepackage{graphicx}
\usepackage{verbatim}
\usepackage{latexsym}
\usepackage{eucal}
\usepackage{tikz}
\usetikzlibrary{matrix}
\usepackage{enumitem}

\numberwithin{equation}{section}

\newtheorem{theorem}{Theorem}[section]
\newtheorem{lemma}[theorem]{Lemma}
\newtheorem{prop}[theorem]{Proposition}

\def \mca {{\mathcal A}}
\def \mcf {{\mathcal F}}
\def \mck {{\mathscr K}}
\def \mcl {{\mathcal L}}
\def \mcr {{\mathcal R}}
\def \mcs {{\mathcal S}}
\def \mct {{\mathscr T}}
\def \mcv {{\mathscr V}}
\def \mbb {{\mathbb B}}
\def \mbc {{\mathbb C}}

\def \mbn {{\mathbb N}}
\def \mbr {{\mathbb R}}

\def \id {\operatorname{Id}}

\def \im {\operatorname{Im}}
\def \re {\operatorname{Re}}
\def \diag{\textrm{Diag}}

\def \beqq {\begin{equation}}
\def \eeqq {\end{equation}}
\def \bpf {\begin{proof}}
\def \epf {\end{proof}}
\def \beq {\begin{equation*}}
\def \eeq {\end{equation*}}

\def \eps {\epsilon}   
\def \la {\lambda}   
\def \La {\Lambda}    
   
\def \lap {\Delta}
\def \p {\partial}
\def \ha {\frac{1}{2}}

\def \tilde {\widetilde}

\def \xo {M\times_0 M}
\def \ff {\text{ff}}
\def \hb {\hbar}
\def \bn {M}
\def \hpl {H_{p_L}}
\def \hpr {H_{p_R}}
\def \bhpl {{}^bH_{p_L}}
\def \bhpr {{}^bH_{p_R}}
\def \rbhpl {\textsf{H}^L_{p}}
\def \rbhpr {\textsf{H}^R_{p}}
\def \pbmo {{}^{b, \ff}T^*(M_0)}
\def \sx {\texttt{x}}
\def \sy {\texttt{y}}

\begin{document}

\title[Sharp resolvent estimates on AHM]{Sharp resolvent estimates on non-positively curved asymptotically hyperbolic manifolds}
\date{\today}
\author{Yiran Wang}
\address{Yiran Wang 
\newline
\indent Department of Mathematics, Emory University}
\email{yiran.wang@emory.edu}

\begin{abstract} 
We study the high energy estimate for the resolvent $R(\la)$ of the Laplacian on non-trapping asymptotically hyperbolic manifolds (AHM). In the literature, polynomial bound of the form $\|R(\la)\| = O(|\la|^{N})$ for $|\la|$ large  and $\la\in \mbc$ in strips where $R(\la)$ is holomorphic was established for some $N > -1$. We prove the optimal bound  $O(|\la|^{-1})$ under the non-positive sectional curvature assumption by taking into account the oscillatory behavior of the Schwartz kernel of the resolvent. 
\end{abstract}

\maketitle

%%%%%%%%%%%%%%%%%%
\section{Introduction}
Let $M$ be an $n+1, n\geq 2$ dimensional smooth manifold with boundary $\p M$ and $g$ be an asymptotically hyperbolic metric on $M$ in the sense of Mazzeo-Melrose \cite{MM}. This means that there is a boundary defining function $\rho$  such that $G = \rho^2 g$ is a smooth Riemannian metric on the closure $\overline M$ of $M$. In this case, sectional curvatures of $(M, g)$ approach $-1$ along curves approaching $\p M.$ In this work, we assume that $(M, g)$ is simply connected and has non-positive sectional curvatures. Then $(M, g)$ is a complete non-compact manifold with no conjugate points and $M$ is diffeomorphic to $\mbb^{n+1} = \{z\in \mbr^{n+1}: |z| < 1\}$. For simplicity, we take $M = \mbb^{n+1}$. 

Consider the (positive) Laplace-Beltrami operator on $(M, g)$, denoted by $\lap_g.$ It is known that $\lap_g$ is an essentially self-adjoint operator on $L^{2}(M) = L^2(M; dg)$. The essential spectrum of $\lap_g$ is $[n^2/4, \infty)$. We consider  $P = \lap_{g} - n^2/4$ and denote the resolvent by   
\beq
R(\la) = (\lap_{g} - \frac{n^2}{4} - \la^2), \quad \la\in \mbc, \im \la \leq 0. 
\eeq
From the spectral theorem, we know that $R(\la)$ is a bounded operator on $L^2(M)$ if $\im\la <<0.$ In \cite{MM}, Mazzeo and Melrose showed that $R(\la)$ can be continued meromorphically to $\mbc\backslash \frac{i}{2}\mbn$ as bounded operators between weighted $L^2$ spaces for fixed $\la$, see also Guillarmou \cite{Gui}. In the last decade, the homomorphic extension of $R(\la)$ and high energy resolvent estimates have drawn lots of attention. For small metric perturbations of the hyperbolic metric near the infinity, which is a class of non-trapping AHM, Melrose, S\'a Barreto and Vasy  showed in \cite{MSV1} that $\rho^aR(\la)\rho^b$ continues holomorphically to strips $\im\la < M, M>0$, provided $|\la| > K(M), b > \im\la, a > \im\la$, and the following estimate  holds
\beqq\label{eqpoly}
\|\rho^aR(\la)\rho^b v\|_{L^2(M)} \leq C|\la|^{-1+ \frac{n}{2}} \|v\|_{L^2(M)}.
\eeqq
Subsequently, the holomorphic extension and the above estimate are proved for non-trapping AHM and the more general conformally compact manifold in S\'a Barreto-Wang \cite{SaWa}. See also Chen-Hassell \cite{CH} for the estimates on the spectrum.  However, these estimates are not optimal. The optimal bound $O(|\la|^{-1})$ is proved by Vasy  \cite{Va} for non-trapping AHM with even metrics at the infinity in the sense of Guillarmou \cite{Gui}. The goal of this paper is to establish the optimal bound for some non-even metrics. Our main result is 
\begin{theorem}\label{main}
Let $(M, g)$ be an $n+1, n\geq 2$ dimensional simply connected asymptotically hyperbolic manifold with non-positive sectional curvatures. Then $\rho^a R(\la)\rho^b$ continues holomorphically to strips $\im\la < M, M>0$, provided $|\la| > K(M), b > \im\la, a > \im\la$ where $K(M)$ is a constant depending on $M.$ Moreover, there exists $C> 0$ such that the following estimate  holds
\beq
\begin{gathered}
\|\rho^aR(\la)\rho^b v\|_{L^2(M)} \leq C|\la|^{-1} \|v\|_{L^2(M)}.
\end{gathered} 
\eeq 
\end{theorem}

We remark that the holomorphic continuation is known from previous work \cite{SaWa}. The improvement is in the order, and this is obtained by exploring the oscillatory behavior of the Schwartz kernel, suggested in \cite{MSV1}. With the non-positive curvature assumption, we show that the distance function behaves similarly to that of a hyperbolic space and a parametrix of the resolvent can be constructed as in \cite{MSV1} using the distance function. When applying the oscillatory integral estimate, we make essential use of the global behavior of the Lagrangian submanifold associated with the geodesic flow, especially the transversality at $\p M$ in a proper sense. We remark that for general non-trapping AHM, a semiclassical parametrix was constructed in \cite{SaWa}. The Lagrangian has a similar global behavior although it cannot be parametrized globally by the distance function.   We expect that Theorem \ref{main} holds as well. 

The paper is organized as follows. We state the structure of the semiclassical parametrix in Section \ref{sec-para}. Then we study the behavior of the distance function and the underlying Lagrangian in Section \ref{sec-dist} and \ref{sec-nonpos}. The parametrix construction follows from these discussions. Finally, we  prove Theorem \ref{main} in Section \ref{sec-res}.

\section*{Acknowledgment}
The author would like to thank Prof.\ Andr\'as Vasy for suggesting the problem and for many stimulating discussions, especially explaining the importance of transversality in the argument. %Most of the work was done when the author was at Stanford. 

%%%%%%%%%%%%%%%%%%
\section{Structure of the parametrix}\label{sec-para}
We translate the high energy estimate to a semiclassical problem. It suffices to consider $\re\la >> 0$ so we let $h = 1/\re \la $ and assume that $h \in (0, 1)$. We write 
\beq
P - \la^2 = h^{-2}(h^2 P - \sigma^2)
\eeq
where $\sigma = \frac{\la}{\re\la} = 1 + i\frac{\im\la}{\re\la}$.  Then we consider the semiclassical operator 
\beq
P(h, \sigma) = h^2(\lap_{g} - \frac{n^2}{4}) - \sigma^2 
\eeq
and take $\sigma \in (1-c, 1+c) \times (-Ch, Ch) \subset \mbc$ for some $c, C > 0.$ We denote the semiclassical resolvent by $R(h, \sigma) = P(h, \sigma)^{-1}$ and observe that $R(\la) = h^2 R(h, \sigma)$. 

We look for an approximation of $R(h, \sigma)$. By the non-positive curvature assumption, there is a unique distance minimizing geodesic between two points $z, z'\in M$ and the distance function $r(z, z')$ on $M\times M$ is smooth away from the diagonal $\diag = \{(z, z) \in M\times M \}$. A natural candidate for the parametrix is the geometric optics ansatz (or WKB solution)
\beq
G(h, \sigma, z, z') = e^{-i\frac{\sigma}{h} r(z, z')}(U_0(\sigma, z, z') + hU_1(\sigma, z, z') + \cdots)
\eeq
The asymptotic behavior  as $z, z'\rightarrow \p M$ can be understood in the $0$-blown-up space introduced in \cite{MM}. To understand the behavior as $h\rightarrow 0$, especially when $\diag$ meets $h = 0$, we work on a blown-up space from  $M\times M\times [0, 1)$ introduced in \cite{MSV1}. We recall these constructions below and we follow the notations there.  

\begin{figure}[htbp]
\centering
\includegraphics[scale=0.65]{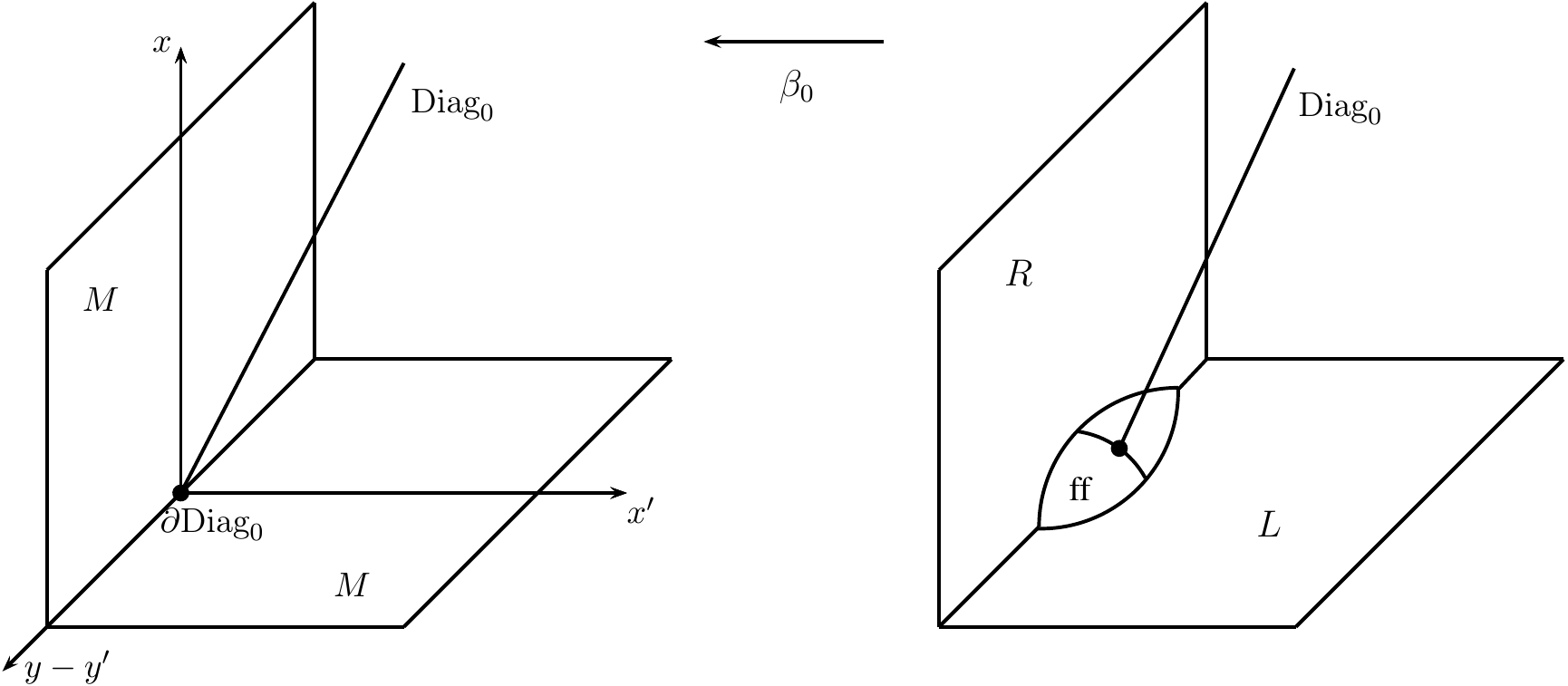}
\caption{The $0$-blown-up space $\bn \times_0 \bn$.}
\label{fig1}
\end{figure}

The first step is to construct the $0$-blown-up of $\bn \times \bn$ as in Mazzeo-Melrose \cite{MM}. Let  $\p \diag = \diag \cap(\p \bn\times \p \bn).$ 
As a set, the $0$-blown-up space (or $0$-stretched product) is 
\beq
M_0 \doteq \bn \times_0 \bn = (\bn\times \bn)\backslash \p\diag \sqcup S_{++}(\p \diag),
\eeq
where $S_{++}(\p\diag)$ denotes the inward pointing spherical bundle of $T_{\p\diag}^*(\bn \times \bn)$. Let 
\begin{gather}
\beta_0: \bn \times_0 \bn \rightarrow \bn \times \bn \label{0blowdown}
\end{gather}
be the blow-down map. Then $\bn \times_0 \bn$ is equipped with a topology and smooth structure  of a manifold with corners for which $\beta_0$ is smooth. The manifold $\xo$  has three boundary hyper-surfaces: the left and right faces $L=\overline{\beta_0^{-1}(\p \bn \times \bn)},$  $R=\overline{\beta_0^{-1}(\bn \times \p \bn)},$  and the front face $\ff= \overline{\beta_0^{-1}(\p\diag)}$. The lifted diagonal is denoted by $\diag_0 = \overline{\beta_0^{-1}(\diag\setminus \p \diag)}$.  It has three co-dimension two corners given by the intersection of  two of these boundary faces, and a co-dimension three corner given by the intersection of  all the three faces. See Figure \ref{fig1}.  
 
\begin{figure}[h]
\centering
\includegraphics[scale=0.65]{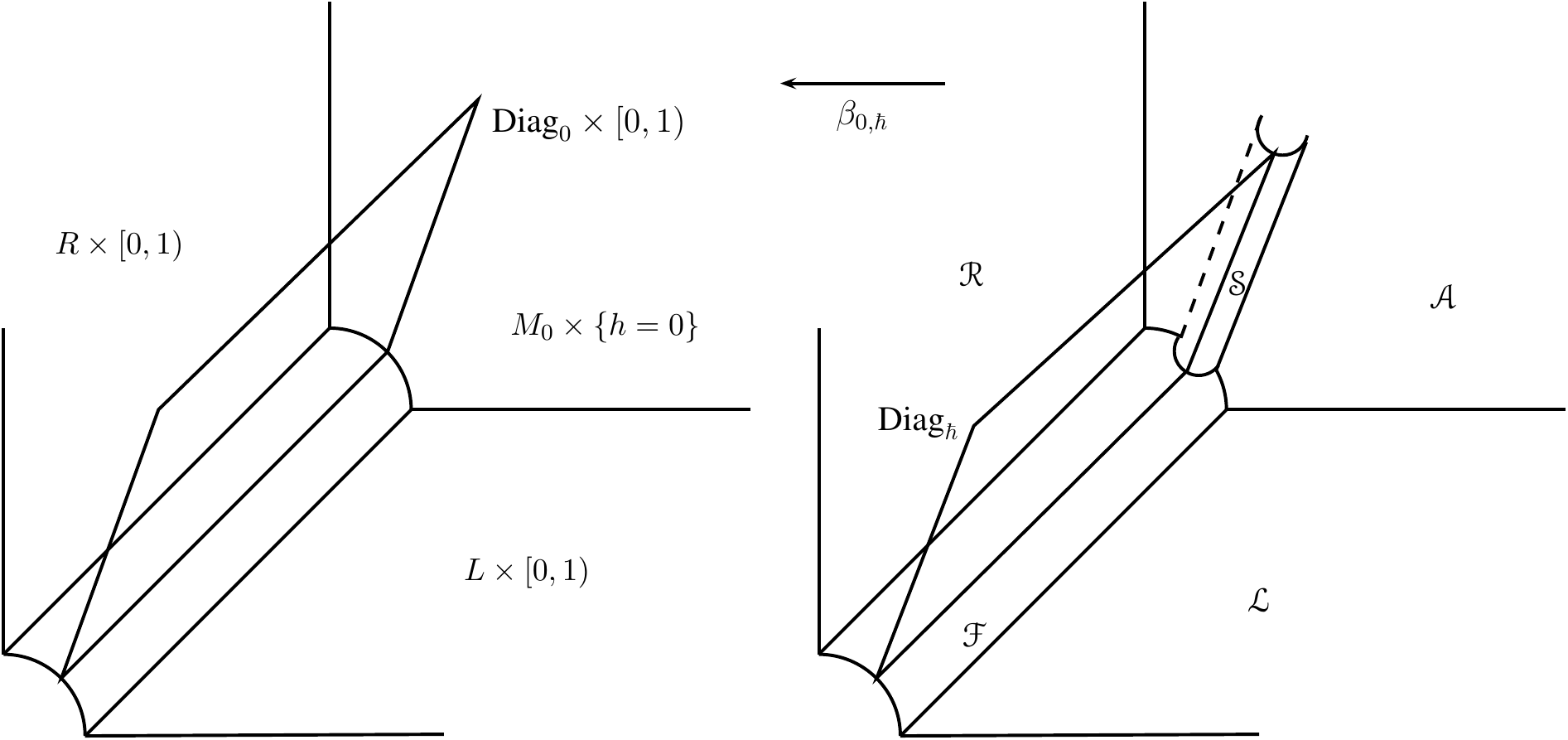}
\caption{The semiclassical blown-up space. The figure on the right is $M_{0, \hbar}$  and the figure on the left is $\bn \times_0 \bn \times [0,1).$}
\label{semiblow}
\end{figure}

Next, on $\xo \times [0, 1)$, the submanifold  $\diag_0 \times[0,1)$  intersect with $\bn \times_0 \bn \times \{0\}$ transversally. We blow up the intersection to get  $M_{0, \hb}$ and denote the associated blow-down map by 
\begin{gather*}
\beta_{0, \hb}: M_{0, \hb}  \longrightarrow \bn \times_0 \bn\times [0, 1).
\end{gather*}
The composition of the blow-down maps $\beta_0$ and $\beta_{0,\hbar}$  will be denoted by
\begin{gather*}
\beta_\hbar = \beta_{0, \hbar}\circ\beta_0 : M_{0, \hbar}  \longrightarrow \bn \times_0 \bn \times [0, 1).
\end{gather*}
Equipped with  a topology so that $\beta_\hbar$ is continuous, the resulting manifold $M_{0, \hbar}$  is a $C^\infty$ compact manifold with corners and it has five boundary faces, see Figure \ref{semiblow}. The left and right faces, denoted by $\mcl, \mcr,$ are the closure of $\beta_{0,\hbar}^{-1} (L \times[0, 1) ), \beta_{0, \hbar}^{-1}(R \times[0, 1))$ respectively.  The front face $\mcf$ is the closure of $\beta_{0, \hbar}^{-1}(\ff\times [0, 1) \backslash (\p \diag_0 \times \{0\})).$ The semiclassical front face $\mcs$ is the closure of $\beta_{0, \hbar}^{-1}(\diag_0\times\{0\})$. Finally, the semiclassical face $\mca$ is the closure of $\beta_{0,\hbar}^{-1}( (\bn \times_0 \bn\backslash \diag_0)\times \{0\})$. The lifted diagonal denoted by $\diag_\hbar$ is the closure of  $\beta_{0, \hbar}^{-1}(\diag_0\times(0, 1))$.  
 
The class of semiclassical pseudo-differential operators on $M_{0, \hbar}$ are defined in two steps. First, let $\Psi_{0, \hbar}^m(\bn)$ consist  of operators $P$ whose kernel $K_P(z, z', h) |dg(z')|$ (as a density) lifts to a conormal distribution of order $m$ to $\diag_\hbar$ and vanishes to infinite order at all faces, except the zero front face $\mcf$, up to which it is smooth with values in conormal distributions and the semiclassical front face $\mcs$, up to which it is $h^{-n-1}C^\infty$. Then we define the space
\beq
\mck^{a, b, c}(M_{0, \hbar}) = \{K \in L^\infty(M_{0, \hbar}): \mcv_b^m K \in \rho_\mcl^a \rho_\mca^b \rho_\mcr^c \rho_\mcs^{-n-1} L^\infty(M_{0, \hbar}), m\in \mbn\},
\eeq
where $\mcv_b$ denotes the Lie algebra of vector fields tangent to $\mcl, \mca, \mcr$. Finally, we define $\Psi_{0, \hbar}^{m, a, b, c}(\bn)$ to be the space consisting of operators $P = P_1 + P_2$ such that $P_1\in \Psi_{0, \hbar}^m(\bn)$ and the  kernel $K_{P_2}|dg(z')|$ is such that $\beta_\hbar^* K_{P_2}\in \mck^{a, b, c}(M_{0, \hbar}).$ 

Now we state the structure of the semiclassical parametrix. 
\begin{prop}\label{thm-para}
Let $(M, g)$ be as in Theorem \ref{main}. Let $r$ be the distance function  of $(M, g)$. 
For $\sigma \in \mbc, \sigma \neq 0, h\in (0, 1)$, there exists a pseudo-differential operator $G(h, \sigma)$ whose kernel is of the form 
\beq
\begin{gathered}
G(h, \sigma, z, z') = G'(h, \sigma, z, z') + G''(h, \sigma, z, z'), 
\end{gathered}
\eeq
where
\beq
\begin{gathered}
G'(h, \sigma) \in \Psi^{-2, \infty, \infty, \infty}_{0, \hbar}(M), \quad G''(h, \sigma, z, z') = e^{-i\frac{\sigma}{h} r(z, z')} U(h, \sigma, z, z')
\end{gathered}
\eeq
with $U(h, \sigma) \in \Psi_{0, \hbar}^{-\infty, \frac n2, -\frac n2-1, \frac n2}(M)$ vanishing in a neighborhood of the semiclassical diagonal $\diag_\hbar$. The operator $G(h, \sigma)$ satisfies  
\beqq\label{eq-paraide}
P(h, \sigma) G(h, \sigma) - \id \in \rho_\mcf^\infty \rho_\mcs^\infty \Psi_{0, \hbar}^{-\infty, \infty, \infty, \frac n2 + i \frac{\sigma}{h}}(M). 
\eeqq
\end{prop}

This result is  the analogue of Theorem 5.1 of \cite{MSV1} which is for small metric perturbations of hyperbolic spaces described as follows. Let $(\mbb^{n+1}, g_0)$ be the Poincar\'e ball model of the hyperbolic space where
\beq
\mbb^{n+1} = \{z\in \mbr^{n+1}: |z| < 1\} \text{ and } g_0 = \frac{4dz^2}{(1 - |z|^2)^2}
\eeq
Consider a one parameter family of perturbations of $g_0$ supported in a neighborhood of $\p \mbb^{n+1}$ of the form 
\beq
g_\delta = g_0 + \chi_{\delta}(z) H(z, dz), 
\eeq
where $H$ is a symmetric $2$-tensor, smooth up to $\p \mbb^{n+1}$ and $\chi_\delta(z) = \chi(\frac{1 - |z|}{\delta})$. Here, $\chi$ is a smooth cut-off function with $\chi(s) = 1, |s|< \ha$ and $\chi(s) = 0, |s|>1$. For $\delta >0$ sufficiently small, the manifold $(\mbb^{n+1}, g_\delta)$ is a non-trapping AHM with no conjugate points. In fact, $x = \frac{1-|z|}{1+|z|}$ is a boundary defining function of $\p \mbb^{n+1}$ and $x^2 g_\delta$ is a non-degenerate Riemannian metric on $\overline{\mbb^{n+1}}.$ 

To prove Proposition \ref{thm-para}, we will show in the next two sections that the behavior of the distance function for a non-positively curved AHM is the same as that for $(M, g_\delta)$. Then the parametrix construction follows the arguments in \cite{MSV1} line by line. The  Lagrangian submanifold associated with the distance function plays an important role when we  apply oscillatory integral estimates later.  

%%%%%%%%%%%%%%%%%%
\section{The underlying Lagrangian submanifold}\label{sec-dist}
Over the interior $M\times M$ and away from $\diag$ where the distance function is smooth,  the graph of  $r$ is a Lagrangian submanifold 
\beq
\La = \{(z, d_z r, z', d_{z'}r): (z, z')\in M\times M\backslash \diag\} \subset T^*(\bn\times \bn \backslash \diag).
\eeq
This can be described as the Hamiltonian flow  as follows. Let $g^*$ be the dual metric of $g$ on $T^*\bn$ and let $p(z, \zeta) =  |\zeta|_{g^*(z)}^2, (z, \zeta)\in T^*M$ be the symbol of $\lap_{g}$. Let $p_L, p_R$ be the lift of $p$ to $\bn\times \bn$ from the left, right factors of $\bn\times \bn$. Let $H_{p_L}, H_{p_R}$ be the associated Hamilton vector fields on $T^*(\bn\times \bn)$. We denote the ``co-sphere" bundle of $\diag$ by
\beq
S^*\diag = \{(z, \zeta, z', -\zeta') \in T^*\bn\times T^*\bn: p(z, \zeta) = p(z', \zeta') = 1\}. 
\eeq 
Then the Lagrangian $\La$ is the joint flow out of $S^*\diag $ of $\hpl, \hpr$
\beq
\begin{gathered}
\La = \bigcup_{t_1, t_2\geq 0} \exp(t_1 \hpl) \circ \exp(t_2 \hpr) (S^*\diag)
\end{gathered}
\eeq
Actually, over the interior $\bn\times \bn$ , one of the flow suffices 
\beq
\La =  \bigcup_{t_1 \geq 0} \exp(t_1 \hpl) (S^*\diag) =  \bigcup_{t_2\geq 0}  \exp(t_2 \hpr) (S^*\diag).
\eeq 

The asymptotic behavior of $\La$ near $\p(M\times M)$ can be understood on the partial $b$-cotangent bundle of $\xo$, denoted by ${}^{b, \ff}T^*(M_0)$, see \cite{MSV1}. It is the dual space of the partially $b$-tangent bundle ${}^{b, \ff}T(M_0)$ whose smooth sections are smooth vector fields on $M_0$ which are tangent to $L, R$ but not necessarily to the front face $\ff$. This can be constructed invariantly via a sequence of blown-ups from the $0$-bundle ${}^0T^*M\times {}^0T^*M$. We refer the readers to \cite{MSV1} for details of the abstract construction and further references on the $b$-geometry and $0$-geometry. The benefit of working with the partial $b$-bundle, as opposed to the $b$-bundle, is that the cotangent bundle over $\ff$ is the usual one and one obtains more precise statement about the regularity of $r$.  

 For our purpose and for readers not familiar with the abstract construction, we find local representations of $\pbmo$. We denote 
\[
\iota_b: T^*(M\times M) \rightarrow {}^{b, \ff}T^*(M_0)
\]
the corresponding bundle map. (This can be regarded as the composition of $\beta_0^*$ and the bundle map $T^*(\xo)\rightarrow \pbmo$.) Let $(x, y_1, \cdots, y_n)$ be local coordinates of $M$ near $\p M$ such that $x$ is a boundary defining function. We use $(x', y')$ for the coordinate of the right factor in $M\times M$.  We denote smooth sections of $T^*(M\times M)$ by $\xi dx + \eta dy + \xi' dx' + \eta' dy'$.  It suffices to find $\iota_b$ over four types of regions in $M_0$ namely (i) near $L$ away from $R, \ff$; (ii) near $L\cap R$ away from $\ff$;  (iii) near $L\cap \ff$ and away from $R$; (iv) near $L\cap R\cap \ff.$

For (i), the vector fields tangent to $L$ is spanned by $x\p_x, \p_{y_j}, \p_{x'}, \p_{y'}$. A basis for $\pbmo$ is $\frac{dx}{x}, dy, dx', dy'$ so that 
\beq
\iota_b:  \xi dx + \eta dy + \xi' dx' + \eta' dy' \rightarrow  \la \frac{dx}{x} + \eta dy + \xi' dx' + \eta' dy', \quad \la = x\xi.
\eeq
(ii) is similar. We have  
\beq
\iota_b:  \xi dx + \eta dy + \xi' dx' + \eta' dy' \rightarrow  \la \frac{dx}{x} + \eta dy + \la' \frac{dx'}{x'} + \eta' dy', \quad \la = x\xi, \la' = x'\xi'.
\eeq
For (iii), we use projective coordinates for the $0$-blown-up. We let $s = x/x', x', y$ and $Y = (y'-y)/x'$ be the local coordinates so $\rho_L = s$ and $\rho_\ff = x'.$ The vector fields tangent to $L$ are spanned by $s\p_s, \p_y, \p_{x'}, \p_{Y'}$ and the dual basis is $\frac{ds}{s}, dy, dx', dY.$ We remark that for the usual $b$-bundle ${}^bT^*M_0$, one would use $\frac{dx'}{x'}$ instead of $dx'$. We compute the image of $\iota_b$ 
\beq
\begin{gathered}
 \xi dx + \eta dy + \xi' dx' + \eta' dy' 
  = \xi(x'ds + sdx') + \eta dy + \xi' dx' + \eta'(Y dx' + x' dY + dy)\\
   = \la \frac{ds}{s} + (\eta + \eta') dy + (\xi' + \eta' \cdot Y )dx' + \eta' x' dY
 \end{gathered}
\eeq
where $\la  = x\xi = sx'\xi.$ 
For (iv), we assume without loss of generality that $y_1' - y_1 >0$ and use projective coordinate near the corner for the $0$-blown-up
\beq
t = y_1' - y_1, s_1 = \frac{x}{y_1'-y_1}, s_2 = \frac{x}{y_1'-y_1}, Z_j=\frac{y_j'-y_j}{y_1'-y_1}, j> 1, y.
\eeq
The boundary surfaces are given by 
\beq
\ff = \{t = 0\}, \quad L = \{s_1 = 0\}, \quad R = \{s_2 = 0\}.
\eeq
The smooth vector fields tangent to $L, R$ are spanned by $s_1\p_{s_1}, s_2\p_{s_2}, \p_t, \p_y, \p_{Z_j}$ so the dual basis is $\frac{ds_1}{s_1}, \frac{ds_2}{s_2}, dt, dy, dZ_j.$ We compute the image of $\iota_b$ as 
\beq
\begin{gathered}
\xi dx + \eta dy + \xi' dx' + \eta' dy'\\
= \xi(s_1 dt + t ds_1) + \eta dy + \xi' (s_2 dt + t ds_2) + \eta'_1(dy_1 + dt) + \sum_{j = 2}^n\eta'_j(Z_j dt + t dZ_j + dy_j)\\
 = \la \frac{ds_1}{s_1} + \la' \frac{ds_2}{s_2} + (\xi s_1 + \xi's_2 + \eta_1 + \eta_j' Z_j)dt + (\eta + \eta') dy +  \sum_{j = 2}^n t \eta_j' dZ_j
\end{gathered}
\eeq
where $\la = x\xi = s_1t\xi, \la' = x'\xi' = s_2t\xi'$.

Let $\bhpl, \bhpr$ be the lift of $\hpl, \hpr$ to $\pbmo$ via $\iota_b$. This can be computed explicitly in local coordinates in case (i)--(iv) just as in \cite[Lemma 2.6]{MSV1}.
% For example, recall that the metric $g$ near $\p M$ can be written as 
%\beq
%g  = \frac{dx^2 + h(x, y, dy)}{x^2} 
%\eeq
%So the symbol $p(x, y, \xi, \eta) = \ha x^2 (\xi^2 + h(x, y, \eta))$. The Hamilton vector field
%\beq
%H_p = \xi x^2 \p_x + x^2 \p_\eta h \p_y - (2x\xi^2 + x^2 h  + x^2 \p_x h) \p_\xi - x^2 \p_y h \p_\eta.
%\eeq
%For (i) and (ii) above, we get $\p_\xi = x\p_\la$ so 
%\beq
%\bhpl = \iota_{b, \ast} (H_p^L) =  x \la \p_x + x^2 \p_\eta h \p_y - (\la^2 +x^2 h(x, y, \eta) + x^2 \p_x h) x \p_\la - x^2 \p_y h \p_\eta.
%\eeq
%The computations near the blow-up are similar to those in \cite{MSV1}.  
%We see that 
%\beq
%\frac{1}{x}\bhpl =   \la \p_x + x  \p_\eta h \p_y - (\la^2 +x^2 h(x, y, \eta) + x^2 \p_x h) \p_\la - x  \p_y h \p_\eta.
%\eeq 
The conclusion (see \cite[Lemma 2.6]{MSV1}) is that these smooth vector fields are tangent to boundary hyper-surfaces of $\pbmo$. More precisely, 
\beq
\bhpl|_{p_L=1} = \pm \textsf{R}_L + \rho_L \tilde W_L, \quad \bhpr|_{p_R = 1} = \pm \textsf{R}_R + \rho_R \tilde W_R
\eeq
where $\textsf{R}_L, \textsf{R}_R$ are the radial vector fields over the left face ${}^bT_L^*(M_0)$, right face  ${}^bT_R^*(M_0)$ and $\tilde W_L, \tilde W_R$ are smooth vector fields tangent to all hyper-surfaces. In particular, the rescaled vector field 
\beq
\rbhpl = \frac{1}{\rho_L}\bhpl, \quad \rbhpr = \frac{1}{\rho_R}\bhpr
\eeq
are transversal to ${}^bT_L^*(M_0)$, ${}^bT_R^*(M_0)$ respectively. The lift of $S^*\diag$ to $\pbmo$ is a submanifold smooth up to $\ff$. 
So up to a re-parametrization,   the lift of the Lagrangian can be described as 
\beq
{}^b\La \doteq \iota_b (\La) =  \bigcup_{t_1, t_2\geq 0} \exp(t_1 \rbhpl) \circ \exp(t_2 \rbhpr) ( \iota_b (S^*\diag))
\eeq
This  is a smooth Lagrangian submanifold outside $\diag_0$ and lies in ${}^bT^*\ff$ over the front face. Here, the symplectic structure on $\pbmo$ is induced from $T^*(M_0)$.

%=======================%
\section{The distance function and parametrix}\label{sec-nonpos}
To construct a parametrix as in \cite{MSV1}, we need to find the asymptotics of the distance function. Let's recall the result for small metric perturbations $(M, g_\delta)$ studied in \cite{MSV1}. Let $r_\delta(z, z')$ be the distance function. For $\delta>0$ sufficiently small and for fixed boundary defining function $\rho_L, \rho_R$ on $M_0$, there is a smooth function $F\in C^\infty(M_0\backslash \diag_0)$ such that 
\beqq\label{eqasym0}
\beta_0^*r_\delta = -\log \rho_L  - \log \rho_R  + F. 
\eeqq
This gives a global parametrization of ${}^b\La$. In \cite{MSV1}, \eqref{eqasym0} is proved by a perturbation argument from the distance function $r_0$ of the standard hyperbolic space $(\mbb^{n+1}, g_0)$ for which the asymptotics can be computed explicitly. We remark that for proper choices of $\rho_L, \rho_R$, $F^2$ is a quadratic defining function of $\diag_0$. If one changes the boundary defining function to $\tilde \rho_L = f \rho_L, \tilde \rho_R = g \rho_R$, where $f, g$ are smooth on $M_0$ and non-vanishing at boundary faces,  we have the same asymptotics but with different $F$
\beq
\beta_0^*r_\delta = -\log \tilde \rho_L  - \log \tilde \rho_R  + \tilde F, \quad \tilde F = F - \log f -\log g. 
\eeq
Here, $\tilde F$ is smooth up to the boundary. In this section, we prove  
\begin{prop}\label{prop-dist}
Let $(M, g)$ be a simply connected AHM with non-positive sectional curvatures. Let $\rho_R, \rho_L$ be boundary defining functions of the right and left face of $\xo$. Then there exists $F\in C^\infty(\xo\backslash \diag_0)$ such that 
\beq
\beta_0^*r = -\log \rho_L  - \log \rho_R  + F. 
\eeq
\end{prop}
 
The problem is near $\p(M\times M)$. We use $(x, y)$ for the local coordinates of $M$ near $\p M$ so that $x$ is a boundary defining function. By a result of Joshi-S\'a Barreto \cite{JS} (see also Graham \cite{Gr}), we can assume  near $\p M$ that 
\beqq\label{eq-metric}
g = \frac{dx^2 + h(x, y, dy)}{x^2}
\eeqq
modulo $O(x^\infty)$ terms, where $h$ is a one parameter family of smooth Riemannian metrics on $\p M$. We first consider the metric in a sufficiently small neighborhood $U$ of a point $(0, y_0)$ at $\p M$. We write 
\beqq\label{eq-h}
h(x, y, dy) = h(0, y_0, dy) + x h_1(x,  y, dy) + \sum_{j=1}^n (y_j - y_{0, j}) h_{2j}(x, y, dy)
\eeqq
where $h_1, h_{2j}$ are smooth on $U.$  
We consider this in $U$ as a small perturbation of the hyperbolic metric 
\beq
\hat g= \frac{dx^2 + h(0, y_0, dy)}{x^2}
\eeq
We know that  $(U, \hat g)$ is isometric to $(U, g_0)$, see \cite{GZ}. Without loss of generality, we identify $\hat g$ with $g_0$ on $U$. Then $g = x^{-2}(x^2 g_0 + \hat h)$ for some two tensor $\hat h$ that can be found from \eqref{eq-h}. For $\delta > 0$, we can shrink $U$ such that $x< \delta, |y_j - y_{0, j}|<\delta$.  The calculation in \cite[Lemma 2.6]{MSV1} can be repeated to show that the Hamilton vector fields  $\rbhpl, \rbhpr$ associated with $g$ are small perturbations of those for $g_0$. Now the perturbation arguments in \cite[Proposition 2.7]{MSV1} apply and we get 
\beq
\beta_0^*r = -\log \rho_L -\log \rho_R + F
\eeq
on the 0-blown-up $U\times_0 U.$ Thus we proved Proposition \ref{prop-dist} near the $0$-blown-up. 
 
 It remains to prove the asymptotics of $r(z, z')$ when $(z, z')$ is away from the diagonal and at least one of $z, z'$ approaches $\p M.$ Now we can work on the product $M\times M$ without worrying about the 0-blown-up. We will use the Jacobi fields and the non-positivity of the curvature to show that the projection of $\La$ to $M\times M$ is a diffeomorphism up to the boundary in the $b$-bundle.

Consider the Hamiltonian flow of $H_p$. This can be identified with the geodesic flow via the bundle isomorphism $T^*M\rightarrow TM$ induced by $g$. We write the dual metric in coordinates $(x, y)$ near $\p M$ as 
\beq
 g^* = x^2 (\p_x^2 + H(x, y, \p_y))
\eeq
Here, $H(x, y, \p_y) = H_{ij}(x, y) \p_{y_i}\p_{y_j}$ and $H_{ij}$ are smooth up to $x = 0.$  
The bundle isomorphism is 
 \beq
\iota: \xi dx +\eta dy \rightarrow x^2 \xi \p_x + x^2 H(x, y)\eta \p_y  
\eeq
or from ${}^bT^*M$ to ${}^bTM$
\beq
\iota: \la \frac{dx}{x} +\eta dy \rightarrow  \la (x\p_x) + x^2 H(x, y) \eta \p_y  
\eeq
In this case, the symbol $p(x, y, \xi, \eta) =  x^2 (\xi^2 + H(x, y, \eta))$ on $T^*M$ and the Hamilton vector field
\beq
H_p = \xi x^2 \p_x + x^2 \p_\eta H \p_y - (2x\xi^2 + x^2 H  + x^2 \p_x H) \p_\xi - x^2 \p_y H \p_\eta.
\eeq
The integral curve $(x(t), y(t), \xi(t), \eta(t)), t\geq 0$ is given by  
\beq
\begin{gathered}
\dot x = \xi x^2, \quad \dot y = x^2 \p_\eta H\\
\dot \xi = - (2x\xi^2 + x^2 H  + x^2 \p_x H) , \quad \dot \eta = - x^2 \p_y H\\
\text{with  $(x(0), y(0), \xi(0), \eta(0)) = (x_0, y_0, \xi_0, \eta_0)$.}
\end{gathered}
\eeq
In particular, $\gamma(t) = (x(t), y(t))$ is a geodesic with $\gamma(0) = (x_0, y_0)\in M, \dot \gamma(0) = \iota(\xi_0, \eta_0) \in T_{(x_0, y_0)}M$. Thus 
\beq
\dot \gamma(t) = \iota(\xi(t), \eta(t)) = (x^2 \xi(t), x^2 H(x, y)\eta(t))   \in TM
\eeq
where $(\xi, \eta)$ satisfies $x^2(\xi^2 + H(x, y, \eta)) = 1.$ So in this setup, $t$ is the arc-length parameter.  One can find how the geodesic equations are related to the system explicitly, see for instance Taylor \cite[Chapter 1, Section 11]{Tay1}. 

Now we look at the geodesics in $b$-bundle. Set $\la  = x\xi$ and the vector field
\beq
\bhpl =   x \la \p_x + x^2 \p_\eta H \p_y - (\la^2 +x^2 H(x, y, \eta) + x^2 \p_x H) x \p_\la - x^2 \p_y H \p_\eta.
\eeq
%We see that 
%\beq
%\frac{1}{x}\bhpl =   \la \p_x + x  \p_\eta h \p_y - (\la^2 +x^2 h(x, y, \eta) + x^2 \p_x h) \p_\la - x  \p_y h \p_\eta.
%\eeq 
The integral curve $(x(t), y(t), \la(t), \eta(t))$ is given by  
\beq
\begin{gathered}
\dot x = x\la, \quad \dot y = x^2 \p_\eta H\\
\dot \la = -x (\la^2 +x^2 H(x, y, \eta) + x^2 \p_x H), \quad \dot \eta = - x^2 \p_y H\\
\text{with  $(x(0), y(0), \la(0), \eta(0)) = (x_0, y_0, \la_0, \eta_0)$}
\end{gathered}
\eeq
and $\la^2 + x^2H(x, y, \eta) = 1.$ For $x = 0,$ we get $\la=\pm 1$. Near $\p M$ where $\la$ is not vanishing, we can use $x$ as the parameter and get 
\beq
\begin{gathered}
\frac{dt}{dx} = \frac{1}{x\la}, \quad \frac{dy}{dx} = \frac{x \p_\eta H}{\la} \\
\frac{d\la}{dx} = - (\la^2 +x^2 H(x, y, \eta) + x^2 \p_x H)/\la, \quad \frac{d\eta}{dx} = - (x  \p_y H)/\la
\end{gathered}
\eeq
In particular, there is a smooth solution $t(x), y(x), \la(x), \eta(x)$ up to $x=0$. From $\frac{dt}{d\ln x} = 1/\la = 1+O(x)$,  we see that $t = \ln x + O(x)$.  So we compactified the arc-length parameter. The integral curve and the geodesic are now extended smoothly to $\p M.$\\
 
Next, we consider the Jacobi field $J$ along a geodesic $\gamma(t), t\geq 0$. In $M$ (the interior), $J$ satisfies the Jacobi equation  
 \beqq\label{eqjac}
 \begin{gathered}
 D_{\dot\gamma(t)}^2 J(t) +R(J(t), \dot \gamma(t)) \dot\gamma(t) = 0\\
 J(0) =V_0, \quad D_{\dot \gamma}J(0) = \dot V_0
 \end{gathered}
 \eeqq
 in which $V_0, \dot V_0 \in TM$, $R$ is the Riemann curvature tensor on $(M, g)$ and $D_{\dot\gamma(t)}$ is the covariant derivative along $\gamma(t)$. As for the geodesics, we look at the equation in ${}^bTM.$ In the follows, for a vector field $V\in TM$ we denote its corresponding $b$-vector field by ${}^b V\in {}^bTM$ via the identification $TM\simeq {}^bTM$, and vice-versa.  
 \begin{lemma}
 The Jacobi field ${}^bJ$ corresponding to $J$ in \eqref{eqjac} has a smooth extension to $\p M$ in ${}^bT\overline M.$
 \end{lemma}
\bpf
 We write the Jacobi equation \eqref{eqjac} as a first order linear system 
 \beqq\label{eqjacnew}
 \begin{gathered}
 D_{\dot\gamma}(J)=  W\\
 D_{\dot\gamma}(W) +R(J, \dot \gamma) \dot\gamma = 0\\
J(0) =V_0, \quad W(0) =  \dot V_0
 \end{gathered}
 \eeqq
where $V_0, \dot V_0\in TM$ and $W\in TM$. Now we rewrite the equation in ${}^bTM$ with  initial data ${}^b V_0, {}^b \dot V_0\in {}^b TM$. 
We show that the solution ${}^b J, {}^b W$ of the system are solutions smooth up to $x = 0$.  

Let's write $\dot \gamma = x^2 \xi \p_x + x^2 H(x, y) \eta \p_y$ where $(\xi, \eta) \in T^*M.$ We have for $k = 0, 1, \cdots, n$ 
 \beqq\label{eqjac1}
 \begin{gathered}
%(\nabla_{\xi \p_x + \eta \p_y}J)^0 = \xi \p_x J^0  + \xi \eta^j \Gamma_{0j}^0 + \eta^i\eta^j \Gamma_{ij}^0\\
%(D_{\dot\gamma(t)}J)^k = (\nabla_{x^2\xi \p_x + x^2\eta \p_y}J)^k = (x^2\xi \dot J^k + x^2\eta^j \dot J^k) + x^2\xi x^2\eta^j \Gamma_{0j}^k + x^2\eta^i x^2\eta^j \Gamma_{ij}^k\\
(D_{\dot\gamma}  J)^k =   \dot J^k +  \sum_{j = 0}^3x^2\xi J^j \Gamma_{0j}^k + x^2\sum_{i = 1}^3 \sum_{j = 0}^3 (H(x, y)\eta)^i  J^j\Gamma_{ij}^k\\
 =  \dot J^k +   \sum_{j = 0}^3 x \la J^j \Gamma_{0j}^k + x^2\sum_{i = 1}^3 \sum_{j = 0}^3 (H(x, y)\eta)^i  J^j\Gamma_{ij}^k
\end{gathered}
 \eeqq
 We examine the coefficients. We already know  that $\la, \eta$ have smooth extensions to $x = 0$ in ${}^bT^*M$. 
The Christoffel symbol is 
 \beq
 \Gamma_{ij}^k = \ha g^{kl}(\p_i g_{lj} + \p_j g_{il} - \p_l g_{ij})
 \eeq
 Because $g_{ij} = O(x^{-2})$ and $g^{ij} = O(x^2)$, we see that $\Gamma_{ij}^k = O(x^{-1})$ in general. Consider 
  \beq
 \Gamma_{0j}^k= \ha g^{kl}(\p_0 g_{lj} + \p_j g_{0l} - \p_l g_{0j})
 \eeq
 When $k = 0$, we have 
   \beq
 \Gamma_{0j}^0= \ha g^{00}(\p_0 g_{0j} + \p_j g_{00} - \p_0 g_{0j}) = \ha g^{00}  \p_j g_{00} 
 \eeq
 Thus for $j\neq 0,$ $\Gamma_{0j}^0 = 0.$ For $j = 0,$ 
we have $\Gamma_{00}^0 =  \ha x^2 \p_x (1/x^2) = -1/x$.  
  Therefore, $ x\la  J^j\Gamma_{0j}^0 = -\la ({}^b J)^0.  $  When $k\neq 0$ but $j = 0$, we get 
    \beq
 \Gamma_{00}^k= \ha g^{kl}(\p_0 g_{l0} + \p_0 g_{0l} - \p_l g_{00}) = 0. 
 \eeq
  
  We look for the equations for ${}^bJ.$  Recall that ${}^bJ = (J^0/x, J^1, J^2, J^3).$ 
  For $k = 0$ in \eqref{eqjacnew} and using \eqref{eqjac1}, the equation for $({}^b J)^0$ is  
 \beqq\label{eq-bJ}
 \begin{gathered}
 %\frac{d(x ({}^b J)^0)}{dx} +  x \la ({}^b J)^j \Gamma_{0j}^0 + x^2(H(x, y)\eta)^i ({}^b  J)^j\Gamma_{ij}^0 = x ({}^b W)^0\\
 x\frac{d ({}^b J)^0}{dx} + ({}^b J)^0 -  \la ({}^b J)^0 + \sum_{i, j = 1}^3 x^2(H(x, y)\eta)^i ({}^b J)^j\Gamma_{ij}^0 = x ({}^b W)^0
 \end{gathered}
 \eeqq
We know that $\la\rightarrow \pm 1$ as $x\rightarrow 0$ and the sign depends on the orientation of the parametrization. Let's consider geodesic approaching $\p M$ as $x\rightarrow 0$. So $\la \rightarrow 1$.  We can divide the equation \eqref{eq-bJ} by $x$ and get 
  \beq
 \begin{gathered}
 \frac{d ({}^b J)^0}{dx}  = A_0({}^b J, {}^b W)
 \end{gathered}
 \eeq
 where $A_0$ is linear in ${}^b J, {}^b W$ with smooth coefficients. For $k\neq 0$, we get from \eqref{eqjacnew} and \eqref{eqjac1} that  
 \beq
  {}^b\dot J^k +  \sum_{j = 1}^3 x \la ({}^b J^j) \Gamma_{0j}^k + x^2\sum_{i =1}^3 \sum_{j = 0}^3 (H(x, y)\eta)^i J^j\Gamma_{ij}^k = ({}^bW)^k
 \eeq
 This can be written as 
   \beq
 \begin{gathered}
 \frac{d ({}^b J)^k}{dx}  = A_k({}^b J, {}^b W)
 \end{gathered}
 \eeq
 where $A_k$ is linear in ${}^b J, {}^b W$ with smooth coefficients. 
 
 Next, for the second equation in \eqref{eqjacnew}, we consider the Riemann curvature tensor 
 \beq
 R_{ij}{}^{k}{}_{l} = \p_i \Gamma_{jl}^k - \p_j\Gamma_{il}^k + \Gamma_{im}^k\Gamma^m_{jl} - \Gamma^k_{jm}\Gamma^m_{il} 
 \eeq
We observe that $ R_{ij}{}^{k}{}_{l}= O(x^{-2})$, hence $R(\cdot, \dot \gamma(t)) \dot\gamma(t) = O(1)$. 
Now the first order system \eqref{eqjacnew} can be written as   
  \beq
 \begin{gathered}
{}^b \dot J  = A({}^b J, {}^b W)\\
 {}^b \dot W =  B({}^b J, {}^b W)
 \end{gathered}
 \eeq
 where $A, B$ are first order linear differential operators with smooth coefficients up to $x = 0.$ Near $\p M$, we use $x$ as the parameter. The geodesic $\gamma$ has a smooth extension to $x = 0$. Hence, the Jacobi field can be smoothly extended to $\p M$ for any given initial data $V_0, \dot V_0$ in ${}^bTM$. This finishes the proof of the lemma. 
\epf
 
 Finally, we consider the geometric implications. Consider the exponential map $\exp: T_pM\rightarrow M, p\in M$ which can be described by $q = \exp(tV) = \gamma(t)$, where $\gamma(t)$ is the unit speed geodesic with $\gamma(0) = p, \dot \gamma(0)= V$ and $|V|_g = 1$. Here, we identified the fiber $T_p M\backslash\{0\}$ with $(0, \infty)\times S_pM$ using the polar coordinate. Now we look at the fiber in ${}^bT_pM$ and consider its compactification. We define
 \[
 q = {}^b \exp(x ({}^bV)) = \gamma(x)
 \]
 where $\gamma(x)$ is the geodesic in ${}^bTM$ and parametrized by $x\in (0, x_V]$. We remark that the parameter $x$ was originally considered near $\p M$. Now we extended it to $(0, x_V]$ along the geodesic $\gamma(t), t\geq 0$ where $x_V > 0$ depends on $V\in S_pM.$  The map ${}^b\exp_p$ can be extended to $[0, x_V]\times {}^bS_pM$ as we have shown. In other words, we found a compactification $\mct_p M$ of ${}^bT_pM$   such that ${}^b\exp_p$ is smooth there. This is also true for $p \in \overline M.$
 \begin{lemma}
 For $p\in \overline M$, ${}^b\exp: \mct_p M\rightarrow \overline M$ is diffeomorphism if and only if the Jacobi field does not vanish. 
 \end{lemma}
 \bpf
 Over the fiber $T_p M, p\in M$, ${}^b \exp$ is the usual exponential map and the result is well known, see \cite[Proposition 10.11]{Lee}. Here, we extend that proof to $\p M.$ By the inverse function theorem, ${}^b \exp$ is a local diffeomorphism near $(x_q, {}^bV) \in [0, x_V]\times {}^bS_pM$ if and only if the push forward ${}^b\exp_*$ is an isomorphism at $(x_q, {}^bV)$. For ${}^bW\in {}^bT_pM$, we have
 \[
 ({}^b\exp_p)_*({}^bW) = \frac{d}{ds}|_{s = 0} {}^b \exp_p(x_q ({}^bV) + s({}^bW))
 \]
 To compute the term, we consider the variation through geodesics $\Gamma_W(s, x) = {}^b\exp_p(x({}^bV + s({}^bW)))$. Over the interior  $M$, the variation is $J_W(x) = \p_s \Gamma_W(0, x)$ a Jacobi field along $\gamma(x)$ and $J_W(x_q) = (\exp_p)_\ast W$. But from Lemma \ref{eqjac}, we know that this can be extended to $\p M$ in the $b$ sense. Thus the map ${}^b\exp_\ast $ is an isomorphism if and only if $J_W(x_q) \neq 0.$ 
 \epf
 
At last, we use the curvature assumption to show that the Jacobi field does not vanish along geodesics approaching $\p M$. Over any compact region of $M$, the curvature is non-negative so the Jacobi fields do not vanish. It suffices to consider near $\p M$ and assume that the curvature is bounded from above by $\delta = -1/R^2$ with $R>0$ close to $1$. We recall the standard proof for conjugate point comparison theorem, see  \cite[Theorem 11.2]{Lee}. For a unit speed geodesic $\gamma(t)$ with $\gamma(0)\in M$,  
 \beq
 \frac{d^2}{dt^2}|J|_g \geq -\frac{R(J, \dot \gamma, \dot \gamma, J)}{|J|_g}\geq - \delta |J|_g
 \eeq
 provided $J\in TM$ and $|J|_g\neq 0$.  We obtain that for $|J(0)|_g =0$
 \beq
 |J|_g \geq R\sinh (t/R) |D_{\dot \gamma}J(0)|_g \geq C R\sinh (t/R) 
 \eeq
 assuming $D_{\dot \gamma}J(0)\neq 0$. We remark that although we showed that $J$ has a smooth extension to ${}^bT\overline M$, $|J|_g$ is not necessarily finite. Consider the geodesic $\gamma(x)$. We re-parametrize it using $\tilde x = x^{1/R}$. This is not smooth at $x =0$ but this is irrelevant in the current argument. We get $t = -R\ln (\tilde x) + O(\tilde x^R)$ and along the geodesic $\gamma(\tilde x)$, we have
\beq
\tilde x |J|_g =  ((\tilde x({}^bJ)^0)^2 + h(\tilde x, y, J^i))^\ha \geq C > 0 
\eeq
which implies that in $b$-bundle, $J^i, i = 1, \cdots, n$ hence ${}^bJ$ is non-vanishing. Since the geodesics are smoothly extended to ${}^bTM$, we can solve the Jacobi equation from the boundary of ${}^bTM.$ Thus the above analysis holds for Jacobi fields along geodesics from $\p M$ to $\p M.$ This shows that the exponential map ${}^b\exp: \mct_p \overline M\rightarrow \overline M$ is a diffeomorphism. \\

Now we conclude that the Lagrangian ${}^b\La$ is diffeomorphic to $M\times M$ away from $\diag $. Over the interior $M\times M\backslash \diag$, this is parametrized by the distance function and near $\p(M\times M\backslash \diag)$, we have
\beq
{}^b\La = \{(x, y, x\p_x r, \p_y r; x', y', x'\p_{x'}r, \p_{y'}r): r = r(x, y, x', y'), (x, y, x', y')\in M\times M\backslash \diag\}
\eeq
But this can be extended to $\p(\xo)$. Hence it follows from the regularity of ${}^b\La$ in ${}^bT^*(\xo)$ that 
\beq
\beta_0^*r = -\log \rho_L -\log \rho_R + F, \quad F \in C^\infty (\xo\backslash \diag_0).
\eeq
This completes the proof of Proposition \ref{prop-dist}. 

\bpf[Proof of Proposition \ref{thm-para}]
With Proposition \ref{prop-dist}, we can repeat the construction of parametrix line by line in \cite[Theorem 5.1]{MSV1}.
\epf

%%%%%%%%%%%%%%%%%%
\section{The resolvent estimate}\label{sec-res}
We prove Theorem \ref{main} in this section. First of all, using  Proposition \ref{thm-para} and applying $R(h, \sigma)$ to \eqref{eq-paraide}, we get $
G(h, \sigma) = R(h, \sigma)( \id + E(h, \sigma) ) $  
where $E$ is the remainder term in Proposition \ref{thm-para}. Let $\rho$ be a boundary defining function. For $a, b \in \mbr$, we consider 
\beq
\rho^aG(h, \sigma)\rho^b = \rho^a R(h, \sigma)\rho^b \rho^{-b}( \id + E(h, \sigma) ) \rho^b =  \rho^a R(h, \sigma)\rho^b  ( \id + \rho^{-b} E(h, \sigma)\rho^b)
\eeq
The $L^2$ estimate of $E(h, \sigma)$ can be obtained by using Schur's lemma, see \cite[Proposition 6.3]{MSV1}. For $\frac{\im\sigma}{h} < b < 2 - \frac{\im\sigma}{h}$, we have
\beq
\| \rho^{-b} E(h, \sigma)\rho^b\|_{L^2} \leq Ch \|f\|_{L^2}.
\eeq
Thus for $h < h_0$ sufficiently small, $\id + \rho^{-b} E(h, \sigma)\rho^b$ is invertible on  $L^2(M)$ so we get 
\beq
\rho^a R(h, \sigma)\rho^b = \rho^aG(h, \sigma)\rho^b ( \id + \rho^{-b} E(h, \sigma)\rho^b)^{-1}
\eeq
It remains to show that 
\beqq\label{eqest}
\|\rho^aG(h, \sigma)\rho^bf\|_{L^2} \leq C h^{-1}\|f\|_{L^2} 
\eeqq
for proper $a, b >\im\sigma/h.$ 
After this is done, we get estimates of $\rho^a R(\la) \rho^b$ via $R(\la) = h^2 R(h, \sigma).$ 
 
We prove \eqref{eqest} using the oscillatory nature of the kernel. Away from $\diag_\hbar$ and the boundary faces, the parametrix in Proposition \ref{thm-para}  is an oscillatory integral. Its prototype on Euclidean space $\mbr^n$ is 
\beq
T_\la(f) = \int_{\mbr^n} e^{i\la \phi(x,y)}a(x, y) f(y)dy, \quad x, y\in \mbr^n
\eeq
where $\phi$ is a real  valued $C^\infty$ (non-homogeneous) phase function and $a$ is smooth with compact support.   
The standard oscillatory integral estimates, see e.g.\ Sogge \cite[Theorem 2.1.1]{Sog} says that if $\phi$ satisfies the non-degeneracy condition 
\beqq\label{eqnondeg}
\det (\frac{\p^2 \phi}{\p x_j \p y_k}) \neq 0
\eeqq
on the support of $a(x, y)$, then for $\la >0$, we have
\beq
\|T_\la f\|_{L^2(\mbr^n)} \leq C \la^{-\frac n2}\|f\|_{L^2(\mbr^n)}. 
\eeq
The non-degeneracy condition can be interpreted geometrically using the associated non-homogeneous Lagrangian submanifold, see \cite[Section 2.1]{Sog}. Let 
\beq
\La_\phi = \{(x, \phi_x'(x, y), y,  \phi_y'(x, y))\} \subset T^*\mbr^n \times T^*\mbr^n
\eeq
be the canonical relation associated to the non-homogeneous phase function. This is a Lagrangian submanifold with respect to the symplectic form $d\xi \wedge dx + d\eta \wedge dy$ on $T^*(\mbr^n\times \mbr^n).$ Consider the two projections
\begin{center}
\begin{tikzpicture}
  \matrix (m) [matrix of math nodes, row sep=1.5em, column sep=1em, minimum width=1.5em]
  {
     & \La_\phi & \\
    (x, \phi_x'(x, y))\in T^*\mbr^n &   & (y,  \phi_y'(x, y)) \in T^*\mbr^n\\};
  \path[-stealth]
    (m-1-2) edge node [above] {$\pi_1$} (m-2-1)
    (m-1-2) edge node [above] {$\pi_2$} (m-2-3);
\end{tikzpicture}
\end{center}
Then \eqref{eqnondeg} is equivalent to that $\pi_1, \pi_2$ are local diffeomorphisms. 
For the estimate of $G(h, \sigma)$, we will make essential use of the flow-out nature of $\La$ and choose submanifolds that are transversal to the flow for performing the oscillatory integral estimates, as  in H\"ormander \cite[Section 4.3]{Ho}.

We write $G(h, \sigma) = G'(h, \sigma) + G''(h, \sigma)$ as in Proposition \ref{thm-para}. We can arrange that the kernel of $G'(h, \sigma)$ is supported close to $\diag_{\hbar}$. Because $G'(h, \sigma) \in \Psi_{0, \hbar}^{-2}(M)$, we get  for  all $a, b$ that 
\[
\|\rho^a G'(h, \sigma) \rho^b f\|_{L^2(M)} \leq C\|f\|_{L^2(M)}.
\]
Hereafter, $C$ denotes a generic constant. Next consider the kernel of $G''$
\[
\rho^aG''(h, \sigma, z, z')\rho^b =  e^{-i\frac{\sigma}{h} r(z, z')} \rho^aU(h, \sigma, z, z') \rho^b
\]
which is supported away from $\diag_\hbar$ in $M_{0, \hbar}$.  It suffices to obtain the estimates near $\mca$ because otherwise $\rho^a G''(h, \sigma)\rho^b$ is bounded in $L^2(M)$ independent of $h$ by Schur's lemma. We divide the estimate into seven regions where the asymptotic behaviors of the kernel are different. 
 \begin{enumerate}%[label=(\roman*)]
 \item near $\mca$ and away from the boundary faces $\mcl, \mcr, \mcf, \mcs$.
 \item near $\mcl\cap \mca$ and away from $\mcr, \mcf, \mcs$.
 \item near $\mcl\cap \mcr\cap \mca$ and away from $\mcf, \mcs$.
 \item near $\mcl\cap \mcf  \cap \mca$ and away from $\mcs, \mcr$.
 \item near $\mcl\cap \mcr\cap \mcf\cap \mca$ and away from $\mcs$.
 \item near $\mcs\cap \mca$ away from $\mcf$ (which is automatically away from $\mcl, \mcr$)
 \item near $\mcs\cap \mcf\cap \mca$ and away from $\mcl, \mcr$.
 \end{enumerate}
Notice that we should also consider the regions obtained by switching $\mcl$ and $\mcr$ in (2) and (4), however, the analysis are identical. After we establish the estimate  \eqref{eqest} for $G''$  in the above types of regions, the proof of Theorem \ref{main} is completed by the  compactness of $M_{0, \hbar}$.

 %===========%
\subsubsection*{{\bf Region (1):}} This is away from $\mcl, \mcr, \mcf$ so the $\rho^a, \rho^b$ factor does not play a role. Using Proposition \ref{thm-para}, we can write $\rho^aG''(h, \sigma)\rho^b$ as
 \beq
 G_1(h, \sigma, z, z') = e^{i\frac{\sigma}{h} r(z, z')} h^{-\frac n2 - 1} A(z, z', h)
 \eeq
 where $A \in C^\infty((M\times M\backslash \diag)\times [0, 1))$ with compact support.  
 For any $(z_0, z_0')\in \bn\times \bn\backslash \diag$, we choose local coordinates $z = (\sx, \sy), z' = (\sx', \sy')$ such that if we denote the hyper-surfaces $S = \{\sx = 0\}\times M, S' = M\times \{\sx' = 0\}$, then $\hpl$ is transversal to $T^*_S(M\times M)$  and $\hpr$ transversal to $T_{S'}^*(M\times M)$. Points on the Lagrangian $\La$ near $(z_0, z_0')$ can be written as
 \beq
(\p_\sx r)  d\sx + (\p_\sy r)  d\sy + (\p_{\sx'} r) d\sx' + (\p_{\sy'} r) d\sy'.
\eeq
Let $\tilde \La$  be the projection of $\La$ to $T^*S\times T^*S'$ or explicitly  
 \beq
  \tilde \La  = \{(\sy, \p_\sy r(0, \sy, 0, \sy'),  \sy',  \p_{\sy'} r(0, \sy, 0, \sy')): (0, \sy, 0, \sy') \text{ is near } (z_0, z_0')\}.
 \eeq
 Because of the transversality of the flow, $\La$ is the joint flow out of $\tilde \La$. 
For example, we can solve $\p_\sx r$ from $p(\p_\sx r, \p_\sy r) = 1$ and integrate along $\textsf{H}_{p}^L$ (and the same for $\textsf{H}_{p}^R$) to get the Lagrangian. This implies that the projections
 \begin{center}
\begin{tikzpicture}
  \matrix (m) [matrix of math nodes,row sep=1.5em,column sep=1em,minimum width=1.5em]
  {
     & \tilde \La & \\
    (\sy, d_\sy r(0, \sy, 0, \sy')) \in T^*\mbr^n &   & (\sy', d_{\sy'}r(0, \sy, 0, \sy'))\in T^*\mbr^n \\};
  \path[-stealth]
    (m-1-2) edge node [above] {$\pi_1$} (m-2-1)
    (m-1-2) edge node [above] {$\pi_2$} (m-2-3);
\end{tikzpicture}
\end{center}
are local diffeomorphism. In a neighborhood $U$ of $(z_0, z_0')$, the projections $\pi_1, \pi_2$ are still diffeomorphisms. Thus for fixed $\sx, \sx'$, the phase function $r$ is non-degenerate in $\sy, \sy'$ variables. Without loss of generality, we assume $a$ is compactly supported in $U$ which is of the product type: $U = U_\sx \times U_\sy \times U_{\sx'}\times U_{\sy'}$ where $U_\bullet$ are relatively compact sets. Then we have
 \beq
 \begin{split}
 \|G_1(h, \sigma)f(z)\|_{L^2} & =  \bigg(\int_{\bn} |\int_{\bn}e^{i\frac{\sigma}{h} r(z, z')} h^{-\frac n2 - 1} A(z, z', h) f(z') dg(z')|^2 dg(z) \bigg)^\ha\\
 &   \leq h^{-\frac n2 -1} \bigg(\int_{U_\sx} \int_{U_\sy} |\int_{U_{\sx'}} \int_{U_{\sy'}} e^{i\frac{\sigma}{h} r(z, z')}  A(z, z', h) f(z') d\sx' d\sy'|^2 d\sx d\sy \bigg)^\ha
 \end{split}
 \eeq
 where we have absorbed the (smooth) density factor $H$ in $dg = H(\sx, \sy) d\sx d\sy$ to the amplitude function. By Minkowski inequality and oscillatory integral estimate, we obtain that 
   \beq
 \begin{split}
 \|G_1(h, \sigma)f(z)\|_{L^2}  &\leq  h^{-\frac n2 -1} \int_{U_{\sx'}} \bigg( \int_{U_\sx} \int_{U_\sy}  |\int_{U_{\sy'}} e^{i\frac{\sigma}{h} r(z, z')}   A(z, z', h) f(\sx', \sy') d\sy'|^2 d\sy d\sx\bigg)^\ha d\sx' \\
& \leq C h^{-\frac n2 -1} \int_{U_{\sx'}} \bigg( \int_{U_\sx} h^{n}  (\int_{U_{\sy'}} |f(\sx', \sy')|^2 d\sy') d\sx\bigg)^\ha d\sx' 
 \leq C h^{-1} \|f\|_{L^2}.
 \end{split}
 \eeq
 
  %===========%
 \subsubsection*{{\bf Region (2) and (3):} }
For region (2), let $(z_0, z_0')\in \p M\times M$. We use local coordinates $z = (x, y)$ near $\p M$ where $x$ is a boundary defining function of $\p\bn$, and $z'$ for $M$. 
So $\rho_\mcl = x, \rho_\mca = h.$ We can write the kernel of $\rho^a G''(h, \sigma) \rho^b$ as 
\beq
 G_2(h, \sigma, x, y, z') =  h^{-\frac n2 - 1}  e^{i\frac{\sigma}{h} r(x, y, z')} x^{a+ \frac n2} A(x, y, z', h)
\eeq
 where $A$ is smooth up to $x = 0$ and supported near $x = 0.$ Now we use the asymptotics in Proposition \ref{prop-dist}
 \[
\beta_0^*r = -\log x + F(x, y, z')
 \]
 where $F$ is smooth up to $x = 0.$ So we have 
 \beq
 \begin{split}
 G_2(h, \sigma, x, y, z') &=  h^{-\frac n2 - 1}  e^{i\frac{\sigma}{h}(-\log x + F(x, y, z'))}A(x, y, z', h)\\
 &= h^{-\frac n2 - 1}  e^{i\frac{\re\sigma}{h} F(x, y, z')} x^{i\frac{\re\sigma}{h} - \frac{\im\sigma}{h} + a+ \frac n2} e^{-\frac{\im\sigma}{h}F(x, y, z')}A(x, y, z', h)
\end{split}
\eeq
We remark that because $|\im\sigma/h| \leq C$ by our assumption and the oscillatory integral estimates do not involve  derivatives in $h$, we can absorb $e^{-\frac{\im\sigma}{h}F}$ to the amplitude function. Also, $|x^{i\frac{\re\sigma}{h}}| =1$ so we can ignore this factor as well. The remaining is an oscillatory integral with phase $F(x, y, z'). $

We look at the Lagrangian $\La$ in $\pbmo$. Over region (2), points in $\La$ can be expressed by 
\beq
\begin{gathered}
(x\p_xr) \frac{dx}{x} + (\p_y r)dy + (\p_{z'}r)dz' 
 = (-1 + xF_x) \frac{dx}{x} + F_y dy + F_{z'} dz'
\end{gathered}
\eeq
We know that the lifted vector field $\textsf{H}_p^L$ is transversal to ${}^bT_L^*(M_0)$. Now we choose local coordinates $(\sx', \sy')$ near $z'_0$ so that $\hpr$ is transversal to $M\times \{\sx' = 0\}$. Let $\tilde \La$ be the projection of $\La$ to $T^*L\times T^*(M\times \{\sx' = 0\})$. Then we conclude that the projections 
 \begin{center}
\begin{tikzpicture}
  \matrix (m) [matrix of math nodes,row sep=1.5em,column sep=1em,minimum width=1.5em]
  {
     & \tilde \La & \\
    (y, F_y(0, y, 0, \sy')) \in T^*\mbr^n&   & (\sy', F_{\sy'}(0, y, 0, \sy')) \in T^*\mbr^n \\};
  \path[-stealth]
    (m-1-2) edge node [above] {$\pi_1$} (m-2-1)
    (m-1-2) edge node [above] {$\pi_2$} (m-2-3);
\end{tikzpicture}
\end{center}
are bijective. Therefore, the Hessian $F_{y\sy'}$ is non-degenerate at $(z_0, z_0')$ in $(y, \sy')$ variables so is in a neighborhood $U$ of $(z_0, z_0')$. Again, we take $U$ to be the product type.  Now we have
 \beq
 \|G_2(h, \sigma) f(z)\|_{L^2} =  \bigg(\int_{\bn} |\int_{\bn}h^{-\frac n2 - 1}  x^{- \frac{\im\sigma}{h} + a + \frac n2} e^{-\frac{\im\sigma}{h}F(x, y, z')} \tilde A(x, y, z', h)  f(z') dg(z')|^2 dg(z) \bigg)^\ha
 \eeq
 Recall the metric $g$ near $\p M$, see \eqref{eq-metric}. The volume form $dg = H(x, y) x^{-(n+1)} dxdy$ with $H$ smooth up to $x = 0$. We have 
 \beq
 \begin{gathered}
 \|G_2(h, \sigma)f(z)\|_{L^2} \\
 = h^{-\frac n2 -1} \bigg(\int_{U_x} \int_{U_y} |\int_{U_\sx'} \int_{U_\sy'} x^{- \frac{\im\sigma}{h} + a + \frac n2} e^{-\frac{\im\sigma}{h}F(x, y, \sx', \sy')} \tilde A(x, y, \sx', \sy', h)   f(\sx', \sy') d\sx' d\sy'|^2 \frac{dx dy}{x^{n+1}}\bigg)^\ha
 \end{gathered}
 \eeq
where the density factor $H$ is absorbed to $\tilde A$. By Minkowski inequality and oscillatory integral  estimates, we get 
  \beq
  \begin{gathered}
 \|G_2(h, \sigma)f(z)\|_{L^2} \leq C h^{-1} \int_{U_\sx'} \bigg( \int_{U_x} x^{ -2\frac{\im\sigma}{h} +2a + n - n-1} \|f(\sx', \cdot)\|^2_{L^2} dx\bigg)^\ha d\sx'  
\leq Ch^{-1}\|f\|_{L^2}
 \end{gathered}
 \eeq
 if the power of $x$ is less than $1$ so $\frac{\im\sigma}{h} < a.$  

The treatment for region (3) is similar. For $(z_0, z_0')\in \mcl\cap \mcr$, we use $z = (x, y), z' = (x', y')$ as local coordinate so that $\rho_\mcl = x, \rho_\mca = h, \rho_\mcr = x'$. We can write $x^a G''(h, \sigma) x^b$ as 
\beq
G_3(h, \sigma, x, y, z') =  h^{-\frac n2 - 1}  e^{i\frac{\sigma}{h} r(x, y, z')} x^{a+ \frac n2} A(x, y, z', h) (x')^{b+ \frac n2}
\eeq
 Now we use $ r = -\log x -\log x' + F(x, y, z')$ from Proposition \ref{prop-dist} to get 
 \beq
 \begin{gathered}
G_3(h, \sigma, x, y, z') 
 = h^{-\frac n2 - 1}  x^{ - \frac{\im\sigma}{h} + a+ \frac n2}  (x')^{ - \frac{\im\sigma}{h} + b+ \frac n2}e^{-\frac{\im\sigma}{h}F(x, y, z')}\tilde A(x, y, x', y', h)
\end{gathered}
\eeq
The Lagrangian $\La$ in $\pbmo$ consists of  
\beq
\begin{gathered}
(x\p_xr)\frac{dx}{x} + (\p_y r) dy + (x'\p_{x'}r)\frac{dx'}{x'} + (\p_{y'} r)dy' \\
 = (-1 + xF_x) \frac{dx}{x} + F_y dy +  (-1 + x'F_{x'}) \frac{dx'}{x'} + F_{y'} dy' 
\end{gathered}
\eeq
We know that the lifted vector field $\textsf{H}_p^L$ is transversal to ${}^bT_L^*(M_0)$ and $\textsf{H}_p^R$ is transversal to ${}^bT_R^*(M_0)$. Let $\tilde \La$ be the projection of $\La$ to $T^*L\times T^*R$ then 
 \begin{center}
\begin{tikzpicture}
  \matrix (m) [matrix of math nodes,row sep=1.5em,column sep=1em,minimum width=1.5em]
  {
     & \tilde \La & \\
    (y, F_y)  \in T^*\mbr^n&   & (y', F_{y'})  \in T^*\mbr^n\\};
  \path[-stealth]
    (m-1-2) edge node [above] {$\pi_1$} (m-2-1)
    (m-1-2) edge node [above] {$\pi_2$} (m-2-3);
\end{tikzpicture}
\end{center}
are bijective. Therefore, the mixed Hessian $F_{y y'}$ is non-degenerate at $(z_0, z_0')$ in $y, y'$ variables so is in a neighborhood $U$ of $(z_0, z_0')$. Use the same arguments, we get 
 \beq
 \begin{split}
& \|G_3(h, \sigma)f(z)\|_{L^2}  
 = h^{-\frac n2 -1} \bigg(\int_{U_x} \int_{U_y} |\int_{U_{x'}} \int_{U_{y'}} x^{- \frac{\im\sigma}{h} + a + \frac n2}(x')^{- \frac{\im\sigma}{h} + b + \frac n2}  \\
 &\quad \quad \quad \quad \quad \quad \cdot e^{-\frac{\im\sigma}{h}F(x, y, x', y')} \tilde A(x, y, x', y', h)   f(x', y') \frac{dx' dy'}{(x')^{n+1}}|^2 \frac{dx dy}{x^{n+1}}\bigg)^\ha\\
 &\leq C h^{-1} \int_{U_{x'}} \bigg( \int_{U_x} x^{ -2\frac{\im\sigma}{h} +2a -1}(x')^{ -\frac{\im\sigma}{h} + b + \frac n 2 - n-1} \int_{U_{y'}}|f(x', y')|^2 dy'  dx\bigg)^\ha d x' \\
&\leq C h^{-1} \int_{U_{x'}} \bigg( (x')^{ -\frac{\im\sigma}{h} + b + \frac n 2 - n-1} \int_{U_{y'}}|f(x', y')|^2 dy'  \bigg)^\ha d x' \\
&\leq C  h^{-1} (\int_{U_{x'}}   (x')^{ 2(-\frac{\im\sigma}{h} + b + \frac n 2 - \frac{n+1}{2})} dx')^\ha ( \int_{U_{x'}}\int_{U_{y'}}|f(x', y')|^2\frac{dy'  dx'}{(x')^{n+1}})^\ha \leq Ch^{-1}\|f\|_{L^2}
 \end{split}
 \eeq
 provided $\im\sigma/h < a$ and $\im\sigma/h <  b.$

  %===========%
\subsubsection*{{\bf Region (4):}  }
 In this region, we use projective coordinates $s = x/x',  Y = (y-y')/x', x', y'$ so that $L = \{s = 0\}$ and $\ff = \{x' = 0\}$ and $\rho_L = s, \rho_\ff = x'$. The distance function 
  \[
 \beta_0^*r = -\log s + F(s, Y, x', y')
 \]
 where $F$ is smooth up to $s = 0, x' = 0$.  The kernel of $x^a G''(h, \sigma) x^b$ can be written as 
 \beq
 \begin{split}
G_4(h, \sigma, x, y, z')  
& = h^{-\frac n2 - 1}  e^{i\frac{\sigma}{h}(-\log s + F(s, Y, x', y'))} s^{a + \frac n2} A(s, Y, x', y', h) (x')^{a +b}\\
&  = h^{-\frac n2 - 1}  e^{i\frac{\im\sigma}{h} F(s, Y, x', y') } s^{\frac{\im\sigma}{h} + a + \frac n2} \tilde A(s, Y, x', y',  h) (x')^{a +b}
\end{split}
\eeq
 where as in the previous cases, we absorbed the irrelevant factors to $\tilde A.$ We consider the kernel near a point $p_0 \in \mcl\cap \mcf$ that is $p_0 = (s_0, Y_0, x'_0, y'_0)$ with $s_0 = 0, x'_0 = 0$. Because $F$ and $\tilde A$ are smooth up to $\ff$, we consider an extension of $M_{0, \hbar}$ across $\mcf$ and extend $F, \tilde A$ smoothly across the faces so that $\tilde A$ is compactly supported. 
 
Now consider $\La$ in $\pbmo$ near $\mcl\cap \mcf$. From the construction in Section 3, we see that it is given by 
\beq
\begin{gathered}
(-1 + sF_s) \frac{ds}{s} + F_{y'} dy' + F_{x'} dx' + F_Y dY
 \end{gathered}
\eeq 
(Here, we notice that in $\pbmo$ the bundle over $\ff$ is the ordinary cotangent bundle.) We know that $\textsf{H}^L_p$ is transversal to ${}^bT_L^*M_0$. Then we choose local coordinates $(\sx', \sy')$ near $(x'_0, y_0')$ so that $\textsf{H}^R_p$ is transversal to $T_{\{x'=0\}}^*(M_0)$. Let $\tilde \La$ be the projection of $\La$ to $T^*L\times T^*\{\sx' = 0\}$. The two projections
 \begin{center}
\begin{tikzpicture}
  \matrix (m) [matrix of math nodes,row sep=1.5em,column sep=1em,minimum width=1.5em]
  {
     & \tilde \La  & \\
    (Y, F_Y) \in T^*\mbr^n&   & (\sy', F_{\sy'}) \in T^*\mbr^n  \\};
  \path[-stealth]
    (m-1-2) edge node [above] {$\pi_1$} (m-2-1)
    (m-1-2) edge node [above] {$\pi_2$} (m-2-3);
\end{tikzpicture}
\end{center}
are bijective at $(s_0, Y_0, x_0', y_0')$. Thus the mixed Hessian $ \p^2 F /(\p Y_j \p \sy'_k)$ is non-degenerate in a neighborhood $U = U_s\times U_{Y}\times U_{\sx'}\times U_{\sy'}$ of $p_0$ (for the smooth extension) where $U_\bullet$ are relatively compact.  

Before we apply oscillatory integral estimates, we need to deal with the singular factor in $dg$. We let $L_0^2 = L^2(\bn; dxdy)$ so that 
 \[
 \|x^{-\frac{n+1}{2}}f\|_{L^2_0}  \sim \|f\|_{L^2}
 \]
 Instead of estimating $G_4''(h, \sigma)$ on $L^2$, we consider the estimate from $L^2_0$ to $L^2$, which amounts to a change of measure. For $f\in L^2_0$, we have
  \beq
  \begin{gathered}
 \|G_4(h, \sigma)x^{\frac{n+1}{2}}  f(z)\|_{L^2}\\
 =  \bigg(\int_{\bn} |\int_{\bn}h^{-\frac n2 - 1}  e^{i\frac{\im\sigma}{h} F(s, Y, x', y') } s^{\frac{\im\sigma}{h} + a  + \frac n2} \tilde A(s, Y, x', y',  h) (x')^{a + b + \frac{n+1}{2}}  f(x', y') dx'dy'|^2 \frac{ ds dY}{s^{n+1}}   \bigg)^\ha\\
 \leq h^{-\frac n2 -1} \int_{U_{\sx'}} \bigg( \int_{U_s} \int_{U_Y}  |\int_{U_{\sy'}}   e^{i\frac{\im\sigma}{h} F(s, Y, x', y') } s^{-\frac{\im\sigma}{h} + a  -\ha}\\
 \cdot \tilde A(s, Y, x', y',  h) (x')^{a+ b + \frac{n+1}{2}}  f(x', y') d\sy'|^2 dY ds \bigg)^\ha d\sx'  
 \leq C h^{-1}\|f\|_{L^2_0}
 \end{gathered}
 \eeq
 by oscillatory integral estimates, provided that $\im\sigma/h < a$ and $a+ b+ \frac{n+1}{2}$ is an integer so that we can extend $\tilde A (x')^{a+ b + \frac{n+1}{2}}$ smoothly across the $\mcf$ face. Now for $f\in L^2$, $x^{-\frac{n+1}{2}} f \in L^2_0$ and we get 
   \beq
  \begin{gathered}
 \|G_4(h, \sigma)  f(z)\|_{L^2}  = \|G_4(h, \sigma)x^{\frac{n+1}{2}}  x^{-\frac{n+1}{2}}  f(z)\|_{L^2}  \leq C h^{-1}\|x^{-\frac{n+1}{2}} f\|_{L^2_0} \leq C h^{-1} \|f\|_{L^2}.
 \end{gathered}
 \eeq
 
Finally, for general $a'+b' + \frac{n+1}{2} \geq0$ non-integer, we choose $a\leq a', b\leq b'$ so that $a+ b+ \frac{n+1}{2}$ is an integer. Then 
 \beq
 \begin{gathered}
 \|x^{a'}G''(h, \sigma)x^{b'}f\|_{L^2} =  \|x^{a'-a} x^{a}G''(h, \sigma)x^{b} x^{b'-b}f\|_{L^2}   \leq C \|x^{a}G''(h, \sigma)x^{b} \|_{L^2}\leq C h^{-1}\|f\|_{L^2}.
 \end{gathered}
 \eeq

  %===========%
\subsubsection*{{\bf Region (5):} }
 In this region, assume that $y_1-y_1' \geq 0$ and we use projective coordinate for the $0$-blown-up
\beq
t = y_1 - y_1', \quad s_1 = \frac{x}{y_1-y_1'}, \quad s_2 = \frac{x}{y_1-y_1'}, \quad Z_j=\frac{y_j-y_j'}{y_1-y_1'}, j> 1, \quad y'.
\eeq
The asymptotics of the distance function is 
\beq
\beta_0^* r = -\log s_1 -\log s_2 + F(s_1, s_2, t, y', Z).
\eeq
Then the kernel of $x^a G''(h, \sigma) x^b$ in this coordinate reads
\beq
\begin{gathered}
G_5(h, \sigma) = (ts_1)^a e^{i\frac{\sigma}{h}(-\log s_1 -\log s_2 + F)} h^{-\frac n2 - 1} s_1^{\frac n2} s_2^{\frac n2} A(t, s_1, s_2, y', Z)\\
 = h^{-\frac n2 -1} e^{-i\frac{\sigma}{h}F} s_1^{\frac{\im\sigma}{h} + a + \frac n2} s_2^{\frac{\im\sigma}{h} + b + \frac n2} \tilde A(t, s_1, s_2, y', Z) t^{a+b}
 \end{gathered}
\eeq
In this case, the Lagrangian $\La$ consists of 
\beq
(-1 + s_1 F_{s_1})\frac{ds_1}{s_1} + F_t dt + F_Z dZ + (-1 + s_2 F_{s_2})\frac{ds_2}{s_2} + F_{y'}dy'
\eeq
The vector fields $\textsf{H}_p^L, \textsf{H}_p^R$ are transversal to $T^*_{L}M_0$ and $T^*_{R}M_0$. Let $\tilde \La$ be the projection of $\La$ to $T^*L \times T^*R$ so the projections 
 \begin{center}
\begin{tikzpicture}
  \matrix (m) [matrix of math nodes,row sep=1.5em,column sep=1em,minimum width=1.5em]
  {
     & \tilde \La & \\
    (t, Z, F_t, F_Z) \in T^*\mbr^n&   & (y', F_{y'}) \in T^*\mbr^n  \\};
  \path[-stealth]
    (m-1-2) edge node [above] {$\pi_1$} (m-2-1)
    (m-1-2) edge node [above] {$\pi_2$} (m-2-3);
\end{tikzpicture}
\end{center}
are bijective.  Then in region (4), we first estimate $G_5(h, \sigma)x^{\frac{n+1}{2}}$ from $L^2_0$ to $L^2$: 
   \beq
  \begin{gathered}
 \|G_5(h, \sigma)x^{\frac{n+1}{2}}f(z)\|_{L^2}\\
 =  \bigg(\int_{\bn} |\int_{\bn}h^{-\frac n2 - 1}  e^{-i\frac{\sigma}{h}F} s_1^{\frac{\im\sigma}{h} + a + \frac n2} s_2^{\frac{\im\sigma}{h} + b + \frac n2} \tilde A(t, s_1, s_2, y, Z) t^{a+b + \frac{n+1}{2}}  ds_2 dy'  |^2 \frac{ds_1 dt dZ_j}{s_1^{n+1}} \bigg)^\ha 
 \end{gathered}
 \eeq
 Then we can treat the estimate as in region (3) and get
 \beq
 \|G_5(h, \sigma)f\|_{L^2} \leq Ch^{-1} \|f\|_{L^2}
 \eeq
 for $a, b > \im\sigma /h.$ 
 
    %===========%
\subsubsection*{{\bf Region (6) and (7):}} For Region (6) we recall that the semiclassical front face is obtained from the blow-up of $r= 0$ and $h = 0$. Away from $\diag_\hbar$, we can use projective coordinates 
\beq
\rho_\mcs = r, \quad \rho_\mca = h/r
\eeq
for $\rho_\mca < \eps$ for some $\eps>0$, that is $ h < \eps r$.  According to Proposition \ref{thm-para}, we know that the kernel $x^a G''(h, \sigma) x^b$ can be written as 
\beq
\begin{gathered}
G_6(h, \sigma) = e^{-i\frac{\sigma}{h} r} (\frac{h}{r})^{-\frac n2 -1} r^{-n-1} A(\rho_\mcs, \rho_\mca)  
= e^{-i\frac{\sigma}{h}} h^{-\frac n2 -1} r^{-\frac n2} A(\frac{h}{r}, r) 
\end{gathered}
\eeq
where $A$ is smooth in $h/r, r$ and compactly supported on $h < r.$ In terms of $z, z'$ variables on $M\times M$, $A$ is smooth in $z, z'$. So we can apply oscillatory integral estimates just as in Region (1) to get 
\beq
\|G_6(h, \sigma) f\|_{L^2}\leq Ch^{-1} \|f\|_{L^2}.
\eeq

For region (7)  near $\mcs\cap \mcf\cap \mca$, we first use projective coordinates for the 0-blown-up 
\beq
s = x/x', x', y', Y = (y- y')/x'
\eeq
Then we blow up $s = 1, Y = 0, h = 0$ and use 
\beq
\rho_\mcs = ((s-1)^2 + |Y|^2)^{\ha}, \quad \rho_\mca = h ((s-1)^2 + |Y|^2)^{-\ha}
\eeq
where $|\cdot|$ stands for the Euclidean norm. Here, $\rho_\mcf = x'.$ Therefore, the kernel of $\rho^a G''(h, \sigma) \rho^b $ can be written as 
\beq
\begin{gathered}
G_7(h, \sigma) = e^{-i\frac{\sigma}{h} r} \rho_\mca^{-\frac n2 -1} \rho_\mcs^{-n-1} A(\rho_\mcs, \rho_\mca) \\
= e^{-i\frac{\sigma}{h}r} h^{-\frac n2 -1}   ((s-1)^2 + |Y|^2)^{-\frac n2} A(h ((s-1)^2 + |Y|^2)^{-\ha}, ((s-1)^2 + |Y|^2)^{\ha} ) 
\end{gathered}
\eeq
where the amplitude $A$ is smooth. We can assume that $A$ is compactly supported in $x' < \eps,  h <  \eps ((s-1)^2 + |Y|^2)^{\ha}$ 
and $s > \eps$ (away from $\mcr$) and $(|s - 1|^2 + |Y|^2)^\ha > \eps$ (away from the diagonal). 
So the kernel is again an oscillatory integral of the type considered in region (1). We get
 \beq
\|G_7(h, \sigma) f\|_{L^2}\leq Ch^{-1} \|f\|_{L^2}.
\eeq

%%===============================REFERENCE==========================================%


\begin{thebibliography}{99}
\bibitem{CH}  X. Chen, A. Hassell. {\em  Resolvent and spectral measure on non-trapping asymptotically hyperbolic manifolds I: Resolvent construction at high energy.}  Comm. Partial Differential Equations 41,  no. 3, 515--578. (2016).

\bibitem{Gr} R.  Graham. {\em Volume and area renormalizations for conformally compact Einstein metrics.} Rend. Circ. Mat. Palermo 63  (2000): 31-42.

\bibitem{Gui} C. Guillarmou. {\em Meromorphic properties of the resolvent on asymptotically hyperbolic manifolds.} Duke Mathematical Journal 129.1 (2005): 1-37.

\bibitem{GZ} L. Guillop\'e, M. Zworski. {\em Polynomial bounds on the number of resonances for some complete spaces of constant negative curvature near infinity.} Asymptotic Analysis 11.1 (1995): 1-22.

\bibitem{Ho} L. H\"ormander. {\em Fourier integral operators. I.} Acta Mathematica 127.1 (1971): 79.

\bibitem{JS} M. Joshi, A. S\'a Barreto. {\em Inverse scattering on asymptotically hyperbolic manifolds.} Acta Mathematica 184.1 (2000): 41-86.

\bibitem{Lee} J. Lee. {\em Riemannian manifolds: an introduction to curvature.} Vol. 176. Springer Science \& Business Media, 2006.

\bibitem{MM} R. Mazzeo, R. Melrose. {\em Meromorphic extension of the resolvent on complete spaces with asymptotically constant negative curvature.} Journal of Functional Analysis 75.2 (1987): 260-310.

\bibitem{MSV1} R. Melrose, A. S\'a Barreto, A. Vasy. {\em Analytic continuation and semiclassical resolvent estimates on asymptotically hyperbolic spaces.} Comm. Partial Differential Equations 39.3 (2014): 452-511.

%\bibitem{MSV2} R. Melrose, A. S\'a Barreto, A. Vasy. {\em Asymptotics of solutions of the wave equation on de Sitter-Schwarzschild space.} Communications in Partial Differential Equations 39.3 (2014): 512-529.

\bibitem{SaWa} A. S\'a Barreto, Y. Wang. {\em The semiclassical resolvent on conformally compact manifolds with variable curvature at infinity.} Comm. Partial Differential Equations 41.8 (2016): 1230-1302.

\bibitem{Sog} C. Sogge. {\em Fourier integrals in classical analysis.} Vol. 210. Cambridge University Press, 2017.

\bibitem{Tay1} M. Taylor. {\em Partial differential equations I: Basic theory.} Vol. 115. Springer Science \& Business Media, 2013.

\bibitem{Va} A. Vasy. {\em Microlocal analysis of asymptotically hyperbolic spaces and high-energy resolvent estimates.} Inverse Problems and Applications: Inside Out II 60 (2013): 487.
 
\end{thebibliography}
\end{document}